\documentstyle[11pt]{article}

\newcommand{\ff}

\newcommand{\sindex}{\end{tabular}\hspace{1em}\begin{minipage}[t]{\textwidth}
   \begin{tabbing}}

\newenvironment{proof}{\noindent{\bf Proof }}{\hfill\rule{2mm}{2mm}}
\newtheorem{theorem}{Theorem}[section]
\newtheorem{lemma}{Lemma}[section]
\newtheorem{corollary}{Corollary}[section]
\newtheorem{proposition}{Proposition}[section]

\newtheorem{defintion}{Defintion}[section]
\newtheorem{example}{Example}[section]

\newenvironment{theo}{\begin{theorem}\em}{\rm\end{theorem}}
\newenvironment{cor}{\begin{corollary}\em}{\rm\end{corollary}}

\newenvironment{prop}{\begin{proposition}\em}{\rm\end{proposition}}

\newenvironment{defi}{\begin{defintion}\em}{\rm\end{defintion}}
\newenvironment{exa}{\begin{example}\em}{\rm\end{example}}

\begin{document}
\title{
   Permutations containing and avoiding \\
	certain patterns }
\author{
   Toufik Mansour \\
   Department of Mathematics \\
   Haifa University}
\maketitle

\subsection*{\centering Abstract}
{\em
Let $T_k^m=\{\sigma\in S_k\mid\sigma_1=m\}$. We prove that the number of 
permutations which avoid all patterns in $T_k^m$ equals $(k-2)!(k-1)^{n+1-k}$ 
for $k\leq n$. We then prove that for any $\tau\in T_k^1$ (or any $\tau\in T_k^k$), 
the number of permutations which avoid all patterns in $T_k^1$ (or in $T_k^k$) 
except for $\tau$ and contain $\tau$ exactly once equals $(n+1-k)(k-1)^{n-k}$ 
for $k\leq n$. Finally, for any $\tau\in T_k^m$, $2\leq m\leq k-1$, this 
number equals $(k-1)^{n-k}$ for $k\leq n$. These results generalize recent 
results due to Robertson concerning permutations avoiding $123$-pattern and 
containing $132$-pattern exactly once.
\section{Introduction} 
\normalsize 
In 1990, Herb Wilf asked the following: How many permutations, of 
length $n$, avoid a given pattern $\tau$? By pattern-avoiding we 
mean the following: 
A permutation $\alpha\in S_n$ avoids a permutation $\tau\in S_m$ if 
there is no $1\leq i_1<\dots<i_m\leq n$ such that 
$(\alpha_{i_1},\dots,\alpha_{i_m})$ is order-isomorphic to 
$\tau=(\tau_1,\dots,\tau_m)$. In this case, we write 
$\alpha\in S_n(\tau)$ and call the permutation $\tau$ a 
\emph{pattern}. If a permutation $\alpha\in S_n$ avoids all 
patterns $\tau$ in a set $T$, we say that $\alpha$ avoids $T$ 
and write $\alpha\in S_n(T)$.\\

The first case to be examined was the case of permutations avoiding one 
pattern of length $3$. Knuth ~\cite{knuth} found that $|S_n(132)|=|S_n(123)|=c_n$ 
where $c_n$ is the $n$-th Catalan number given by the formula 
$c_n=\frac{1}{n+1}{{2n}\choose n}$, and it is easy to prove that 
$|S_n(\tau)|=c_n$ for all $\tau\in S_3$. Later then Simion and 
Schmidt ~\cite{sim} found the cardinalities of $|S_n(T)|$ for all 
$T\subseteq S_3$.
 
A permutation $\alpha\in S_n$ contains $\tau\in S_m$ exactly $r$ times 
if there exist exactly $r$ different sequences $1\leq i_1^j<\dots<i_m^j\leq n$, 
$1\leq j\leq r$, such that $(\alpha_{i_1^j},\dots,\alpha_{i_m^j})$ 
is order-isomorphic to $\tau=(\tau_1,\dots,\tau_m)$ for all $j=1,\dots,r$.\\

Noonan shows in ~\cite{noon} that the number of permutations 
containing exactly one 123-pattern  is given by a simple formula 
$\frac{3}{n} {{2n}\choose {n+3}}$. Bona ~\cite{bona} proves that expression 
${{2n-3}\choose {n-3}}$ enumerates those permutations containing exactly one 
$132$-pattern, and this result was extended by Robertson, Wilf and Zeilberger 
~\cite{rwz} to calculate the number of $132$-avoiding permutations that have 
a given number of $123$-patterns. 

Robertson ~\cite{robert} proves that $(n-2)2^{n-3}$ gives the number of 
permutations containing exactly one $132$-pattern and avoiding the 
$123$-pattern; he also shows that the number of permutations containing 
exactly one $123$-pattern and one $132$-pattern equals $(n-3)(n-4)2^{n-5}$.\\

In this note we obtain a generalization of Robertson's result concerning 
the number $132$-avoiding permutations that have one $123$-pattern. 
As a byproduct we get a generalization of the following result due to Simion 
and Schmidt ~\cite{sim}: $|S_n(123,132)|=|S_n(213,231)|=2^{n-1}$.\\

{\bf Acknowledgements}. I am grateful to R. Chapman and to the 
anonymous referee for their helpful comments. Special thanks to my advisor 
A. Vainshtein for his help during the preparation of the paper.
\section{Permutations avoiding $T_k^m$}
Let $T_k^m=\{\sigma\in S_k\mid \sigma_1=m\}$ for all $1\leq m\leq k$, 
$k\in\mathcal{N}$. For example $T_3^1=\{123,132\}$ and 
$T_4^2=\{2134,2143,2314,2341,2413,2431\}$.

As a preparatory step, we calculate the cardinalities of the sets 
$S_n(T_k^m)$ for all $1\leq m\leq k$, $k\in\mathcal{N}$.

\begin{theo}
\label{tkm}
	 $$|S_n(T_k^m)|=(k-2)!\cdot(k-1)^{n+2-k},$$ 
	 for all $m,k,n\in\cal N$ such that $2<k\leq n$, $1\leq m\leq k$.
\end{theo}
\begin{proof} 
  Let $G_n=S_n(T_k^m)$, and define the family of functions $f_h:S_n\rightarrow S_{n+1}$ by:
  \begin{center}
  $[f_h(\beta)]_i=\left\{ \begin{array}{ll} h,		& \mbox{when $i=1$} \\
					   \beta_{i-1},	& \mbox{when $\beta_{i-1}<h$} \\
					   \beta_{i-1}+1, & \mbox{when $\beta_{i-1}\geq h$} \end{array} \right.$
  \end{center}
for every $i=1,...,n+1$, $\beta\in S_n$ and $h=1,...,n+1$.\\

From this we see that:
\begin{center} if $\sigma\in G_n$ then 
     $f_{n+1}(\sigma ),f_n(\sigma ),...,f_{n+m-k+2}(\sigma),f_1(\sigma),...,f_{m-1}(\sigma)\in G_{n+1}$ 
\end{center}
so $(k-1)\cdot |G_n|\leq |G_{n+1}|$ where $n\geq k$.\\

Assume that $(k-1)\cdot |G_n|<|G_{n+1}|$. Then there exists a permutation 
$\alpha\in G_{n+1}$ such that $m\leq\alpha_1\leq n+m-k+1$, so there 
exist $k-1$ positions $1<i_1<\cdots <i_{k-1}\leq n+1$ such that the 
subsequence $\alpha_1,\alpha_{i_1},\ldots,\alpha_{i_{k-1}}$ is order-isomorphic to 
one of the patterns in $T_k^m$, which contradicts the definition of 
$G_{n+1}$. So $(k-1)\cdot |G_n|=|G_{n+1}|$ for $n\geq k$. Besides $|G_k|=(k-1)(k-1)!$ 
(from the definition of this set), hence $|G_n|=(k-2)!\cdot (k-1)^{n+2-k}$ 
for $n\geq k$.
\end{proof}

\begin{exa} As a corollary we get the result of Simion and Schmidt 
~\cite{sim}: $|S_n(T_3^1)|=|S_n(T_3^2)|=2^{n-1}$ for all $n\in\cal N$. 
Other examples: $|S_n(T_4^1)|=2\cdot 3^{n-2}$ and 
$|S_n(T_5^2)|=3\cdot 2^{2n-5}$.
\end{exa}

\begin{cor} 
     $$|S_n(T_k^a\cup T_k^{a+1} \cup \dots\cup T_k^b )|=(k-1)!(k+a-b-1)^{n+1-k},$$ 
     where $1\leq a\leq b\leq k$.
\end{cor}
\begin{proof} 
Let $G_n=S_n(T_k^a\cup T_k^{a+1} \cup\dots\cup T_k^b)$. 
From Theorem \ref{tkm} we get that $\alpha\in G_n$ if and only if 
  $$f_1(\alpha),\dots,f_{a-1}(\alpha),f_{n+b-(k-2)}(\alpha),\dots,f_{n+1}(\alpha) \in G_{n+1}.$$

So $|G_{n+1}|=(k+a-b-1)|G_n|$. Besides $|G_k|=(k+a-b+1)(k-1)!$, hence the 
theorem holds.
\end{proof}
\\

For a general family $T_k^{i_1},\dots,T_k^{i_d}$ the situation is more 
complicated. However, the following recurrence can be proved.

\begin{cor} 
 For all $n\geq 2k+1$, 
     $$|S_n(T_k^{i_1}\cup \dots \cup T_k^{i_d} )|=(k+i_1-i_d-1)|S_{n-1}(T_k^{i_1}\cup\dots\cup T_k^{i_d})|,$$ 
     where $1\leq i_1<i_2<\dots< i_d\leq k$.
\end{cor}

\section{Permutations avoiding $T_k^1\backslash{\{\tau\}}$ and containing $\tau$ exactly once}

Let $b_1<\dots<b_n$; we denote by $S_{\{b_1,\dots,b_n\}}$ the set 
of all permutation of the numbers $b_1,\dots,b_n$; for example, 
$S_{\{1,\dots,n\}}$ is just $S_n$. As above we denote by 
$S_{\{b_1,\dots,b_n\}}(T)$ the set of all permutations in 
$S_{\{b_1,\dots,b_n\}}$ avoiding all the permutations in $T$.

\begin{prop}
   Let $T\subseteq S_{\{c_1,\dots,c_k\}}$. Then 
   there exists $R\subseteq S_k$ such that $|S_n(R)|=|S_{\{c_1,\dots,c_n\}}(T)|$.
\label{setword}
\end{prop}
\begin{proof}
We define a function $f:S_{\{c_1,\dots,c_k\}}\rightarrow S_k$ by
   $$f((c_{i_1},c_{i_2},\dots,c_{i_k}))=(i_1,i_2,\dots,i_k),$$ 
then evidently $|S_{\{c_1,\dots,c_n\}}(\tau)|=|S_n(f(\tau))|$ for
all $\tau\in S_{\{c_1,\dots,c_k\}}$.

Let $T=\{\tau_1,\dots,\tau_l\}$ and $R=\{f(\tau_1),\dots,f(\tau_l)\}$. So 
  $$S_{\{c_1,\dots,c_n\}}(T)=\bigcap_{i=1}^l S_{\{c_1,\dots,c_n\}}(\tau_i),$$
hence by the isomorphism $f$ we have that 
   $|S_{\{c_1,\dots,c_n\}}(T)|=|S_n(R)|$.
\end{proof}
\begin{defi} 
  Let $M_{\tau}^{k,m}=T_k^m\backslash\{\tau\}$, for $\tau\in T_k^m$.
  We denote by $S_n(T_k^m;\tau)$ the set of all permutations in 
  $S_n$ that avoid $M_{\tau}^{k,m}$ and contain $\tau$ exactly once.
\end{defi}

\begin{theo} 
\label{mk1}
    $$|S_n(T_k^1;\tau)|=(n+1-k)\cdot (k-1)^{n-k},$$ 
    for all $k\leq n$, $\tau\in T_k^1$.
\end{theo}
\begin{proof} 
  Let $\alpha\in S_n(T_k^1;\tau)$, and let us consider the possible 
  values of $\alpha_1$:
\begin{enumerate}
  \item	$\alpha_1\geq n-k+2$. Evidently $\alpha\in S_n(T_k^1;\tau)$ 
	if and only if $\alpha\in S_{\{1,\dots,n\}\backslash\{\alpha_1\}}(T_k^1;\tau)$.

  \item	$\alpha_1\leq n-k$. Then there exist $1<i_1<\dots<i_k\leq n$ such 
    	that $(\alpha_1,\alpha_{i_1},\dots,\alpha_{i_k})$ is a 
	permutation of the numbers $n,\dots,n-k+1,\alpha_1$.
	For any choice of $k-1$ positions out of $i_1,\dots,i_k$, the 
	corresponding permutations preceeded by $\alpha_1$ is 
	order-isomorphic to some permutation in $T_k^1$. Since 
	$\alpha$ avoids $M_{\tau}^{k,1}$, it is, in fact, 
	order-isomorphic to $\tau$. We thus get at least $k$ occurrences 
	of $\tau$ in $\alpha$, a contradiction.

  \item	$\alpha_1=n-k+1$. Then there exist $1<i_1<\dots<i_{k-1}\leq n$ such that 
	$\eta =(\alpha_1,\alpha_{i_1},\dots,\alpha_{i_{k-1}})$ is a 
	permutation of the numbers $n,\dots,n-k+1$. As above, we 
	immediately get that $\eta$ is order-isomorphic to $\tau$.
	We denote by $A_n$ the set of all permutations in $S_n(T_k^1;\tau)$ 
	such that $\alpha_1=n-k+1$, and define the family of 
	functions $f_h:A_n\rightarrow S_{n+1}$ by:
	    $$[f_h(\beta)]_i=\left\{ \begin{array}{ll} 1,  & \mbox{when $i=h$} \\
					   \beta_i+1,     & \mbox{when $i<h$} \\
					   \beta_{i-1}+1, & \mbox{when $i>h$} \end{array} \right. ,$$
	for every $i=1,...,n+1$, $\beta\in A_n$ and $h=1,...,n+1$.
	It is easy to see that for all $\beta\in A_n$,
	     $$f_{n+1}(\beta),\dots,f_{n-k+3}(\beta)\in A_{n+1},$$
	hence $(k-1)|A_n|\leq |A_{n+1}|$.

        Now we define another function $g:A_{n+1}\rightarrow S_n$ by:
	\begin{center}
  	    $$[g(\beta)]_i=\left\{ \begin{array}{ll}   
				\beta_i-1,     & \mbox{when $i<h$} \\
	                        \beta_{i+1}-1, & \mbox{when $i+1>h$} \end{array} \right. ,$$
	\end{center}
	where $\beta_h=1$, $i=1,...,n$, $\beta\in A_{n+1}$.

	Observe that $h\geq n-k+3$, since otherwise already 
	$(\beta_h,\beta_{h+1},\dots,\beta_{n+1})$ contains a pattern
	from $T_k^1$, a contradiction. It is easy to see that 
	$g(\beta)\in A_n$ for all $\beta\in A_{n+1}$, hence 
	$|A_{n+1}|\leq (k-1)|A_n|$. So finally, 
	$|A_{n+1}|=(k-1)|A_n|$ and $|A_n|=(k-1)^{n-k}$, since $|A_k|=1$.
   \end{enumerate}
Since the above cases $1,2,3$ are disjoint, and by Theorem \ref{tkm} 
and  Proposition \ref{setword} 
we obtain
	$$|S_n(T_k^1;\tau)|=(k-1)|S_{n-1}(T_k^1;\tau)|+(k-1)^{n-k},$$
hence
	$$|S_n(T_k^1;\tau)|=(n+1-k)\cdot (k-1)^{n-k},$$ 
for all $k\leq n$, $\tau\in T_k^1$.
\end{proof}

\begin{exa}
$|S_n(123;132)|=|S_n(132;123)|=(n-2)2^{n-3}$, which is the  
result of Robertson in ~\cite{robert}.
\end{exa}

\begin{cor} 
   $$|S_n(T_k^k;\tau)|=(n+1-k)\cdot (k-1)^{n-k},$$ 
   for all $k\leq n$, $\tau\in T_k^k$. 
\end{cor}
\begin{proof} 
  Let $\beta$ be a permutation complement to $\tau$.
  By the natural bijection between the set $S_n(T_k^1;\beta)$ and 
  the set $S_n(T_k^k,\tau)$ for all $\tau\in T_k^k$ we have that 
  $S_n(T_k^k;\tau)$ have the same cardinality as $S_n(T_k^1;\beta)$, 
  which is $(n+1-k)\cdot (k-1)^{n-k}$ by Theorem \ref{mk1}.
\end{proof}
\section{Permutations avoiding $T_k^m\backslash{\{\tau\}}$ and containing $\tau$ exactly once, $2\leq m \leq k-1$} 

Now we calculate the cardinalities of the sets 
$S_n(T_k^m;\tau)$ where $2\leq m\leq k-1$, $\tau\in T_k^m$.
\begin{theo} 
   $$|S_n(T_k^m;\tau)|=(k-1)^{n-k},$$ 
   for all $2\leq m<k\leq n$, $\tau\in T_k^m$.
\end{theo}
\begin{proof} 
  Let $G_n=S_n(T_k^m;\tau)$, $\alpha\in G_n$, and let us consider the 
  possible values of $\alpha_1$:
\begin{enumerate}
  \item	Let $\alpha_1\leq m-1$. Evidently $\alpha\in G_n$ if and only if 
	$(\alpha_2,\dots,\alpha_n)\in S_{\{1,\dots,n\}\backslash\{\alpha_1\}}(T_k^m;\tau)$.

  \item $\alpha_1\geq n-k+m+1$. Evidently $\alpha\in G_n$ if and only if 
	$(\alpha_2,\dots,\alpha_n)\in S_{\{1,\dots,n\}\backslash\{\alpha_1\}}(T_k^m;\tau)$.

  \item $m\leq\alpha_1\leq n-k+m$. By definition we have that $|G_k|=1$, so 
	let $n\geq k+1$. If $m+1\leq\alpha_1$ then $\alpha$ contains at 
	least $m\geq 2$ occurrences of a pattern from $T_k^m$, and 
	If $\alpha_1\leq n-k+m-1$ then $\alpha$ conatins at least 
	$k-m+1\geq 2$ occurrences of a pattern from $T_k^m$, 
	a contradiction.
  \end{enumerate}
Since the above cases $1,2,3$ are disjoint and by Proposition 
\ref{setword} we obtain $|G_n|=(k-1)|G_{n-1}|$ for all $k\leq n$. 
Besides $|G_k|=1$, hence $|S_n(T_k^m;\tau)|=(k-1)^{n-k}$.
\end{proof}

\begin{exa} 
$|S_n(213;231)|=|S_n(231;213)|=2^{n-3}$. 
\end{exa}

\end{document}